\documentclass[format=sigplan,9pt,authorversion]{acmart}

\usepackage{tikz}
\usetikzlibrary{arrows}

\pgfmathsetmacro{\arrow}{1.3}   

\tikzset{crosses over/.style={
    preaction={draw, line width=0.2cm, white,-}
  }
}

\newcommand{\poly}{\mathsf{Poly}}
\newcommand{\dill}{\mathsf{Dill}}

\newtheorem*{theorem*}{Theorem}

\begin{document}

\acmConference[LICS]{Logic in Computer Science}{9--12 July 2018}{Oxford, UK}

\title{Dialectica models of type theory}

\author{Sean K. Moss}
\affiliation{%
  \department{Department of Computer Science}
  \institution{University of Oxford}
}
\email{sean.moss@univ.ox.ac.uk}

\author{Tamara von Glehn}
\affiliation{%
  \department{Department of Pure Mathematics and Mathematical Statistics}
  \institution{University of Cambridge}
}
\email{t.l.von-glehn@dpmms.cam.ac.uk}

\begin{abstract}
  We present two Dialectica-like constructions for models of intensional Martin-L\"of type theory based on G\"odel's original Dialectica interpretation and the Diller-Nahm variant, bringing dependent types to categorical proof theory.
  We set both constructions within a logical predicates style theory for display map categories where we show that `quasifibred' versions of dependent products and universes suffice to construct their standard counterparts.
  To support the logic required for dependent products in the first construction, we propose a new semantic notion of finite sum for dependent types, generalizing finitely-complete extensive categories.
  The second avoids extensivity assumptions using biproducts in a Kleisli category for a fibred additive monad.
\end{abstract}

\copyrightyear{2018}
\acmYear{2018}
\setcopyright{acmlicensed}
\acmConference[LICS '18]{33rd Annual ACM/IEEE Symposium on Logic in Computer Science}{July 9--12, 2018}{Oxford, United Kingdom}
\acmBooktitle{LICS '18: 33rd Annual ACM/IEEE Symposium on Logic in Computer Science, July 9--12, 2018, Oxford, United Kingdom}
\acmPrice{15.00}
\acmDOI{10.1145/3209108.3209207}
\acmISBN{978-1-4503-5583-4/18/07}

\maketitle

\section{Introduction}\label{sec:introduction}

G\"odel's `Dialectica interpretation' \cite{GODEL_dialectica1958} is an interpretation of Heyting arithmetic into a system of computable functionals, establishing the relative consistency of the former \cite{TROELSTRA1973,KOHLENBACH2008,AVIGADFEFERMAN1998}.
The most characteristic aspects of the interpretation are its reduction of every proposition to one of the form $\exists \vec u \forall \vec x A(\vec u;\vec x)$ where $A$ is quantifier-free (a `Dialectica proposition'), and its handling of the implication.
For our purposes, it suffices to understand the handling of implication as a method for reducing a formula
\begin{equation}\label{eq:1}
  \Big(
  \exists u^U \forall x^X A(u;x)
  \Big)
  \to
  \Big(
  \exists v^V \forall y^Y B(v;y)
  \Big),
\end{equation}
where $A$ and $B$ are quantifier-free and $u,v,x,y$ are variables with specified types, to a Dialectica proposition $\exists \vec w \forall \vec z\, C(\vec w;\vec z)$ (where $C$ is quantifier-free).
Following Dialectica, the reduction of \eqref{eq:1} would be
\vspace{-0.15cm}
\begin{equation}
  \label{eq:2}
  \exists f^{U \to V}, F^{U \times Y \to X} \forall u,y \Big(A(u;F(u,y)) \to B(f(u);y)\Big).
\end{equation}
To see why this is reasonable, consider that a constructive realization of the implication \eqref{eq:1} should, in particular, map potential witnesses $u$ of the antecedent to potential witnesses $v = f(u)$ of the consequent, and that given such $u$ and $v$ it should also map potential counterexamples $y$ of the consequent to potential counterexamples $x = F(u,y)$ of the antecedent.

Here we present some of the results of the authors' PhD theses \cite{VONGLEHN2014,MOSS2017}, which aim to give Dialectica-style functional interpretations for intensional Martin-L\"of type theory (dependent type theory with at least $\Sigma$-, $\Pi$-, and identity types) \cite{MARTINLOF1984}.
Here we treat both G\"odel's original Dialectica and the variant due to Diller and Nahm \cite{DILLERNAHM1974} with our main results.

\begin{theorem*}[\ref{theorem-polynomial}]
  A type theory with $\Sigma$-, $\Pi$-, identity, and strongly extensive finite sum types and a universe gives rise to a type theory $\poly$ whose types are `Dialectica propositions', also modelling $\Sigma$-, $\Pi$-, identity, and strongly extensive finite sum types and a universe.
\end{theorem*}

\begin{theorem*}[\ref{theorem-diller-nahm}]
  A type theory with $\Sigma$-, $\Pi$-, identity, finite sum, and finite multiset types and a universe gives rise to a type theory $\dill$ whose types are `Diller-Nahm propositions', also modelling $\Sigma$-, $\Pi$-, and identity types and a universe.
\end{theorem*}

By implementing type constructors in the resulting type theory $\poly$ of \ref{theorem-polynomial}, we can define interpretations of other type theories into $\poly$.
These could be used, for example, to give relative consistency results, generalizing G\"odel's original argument to dependent types.
Actually, since the situation is already quite complicated, here we will drop the layer of predicates from our Dialectica (and Diller-Nahm) propositions, which is to say we consider only the vectors of types of the quantified variables.
Equivalently, we only consider those Dialectica propositions of the form $\exists \vec u \forall \vec x\, \top$.
We call the model $\poly$ of Theorem \ref{theorem-diller-nahm} the \emph{polynomial model} because its underlying category is the category of non-indexed \emph{polynomials} or \emph{containers} (see \cite{AltenkirchLevyStaton2010,AbbottAltenkirchGhani2003}).

In \cite{HYLAND2002}, categorical proof theory is propounded as a lens through which to study the structure of proofs and the machinery of proof theory via the insights of category theory.
With this work we continue the strand of investigation launched by de Paiva \cite{DEPAIVA_thesis1991}, who gave the construction of a \emph{Dialectica category} whose morphisms correspond to realizations of \eqref{eq:2} in some category of types.
Under certain assumptions on the basic category $\mathbb C$, the Dialectica category $\mathsf{Dial}(\mathbb C)$ is symmetric monoidal closed, and hence a model for (propositional) linear logic \cite{GIRARD_linearlogic1987}.
Moreover, under further assumptions, $\mathsf{Dial}(\mathbb C)$ together with a certain comonad is a model of linear logic with the `$!$' modality.
The Kleisli category $\mathsf{Dial}(\mathbb C)_!$ for this comonad is a `Dialectica category' but for the Diller-Nahm variant rather than G\"odel's original interpretation, and turns out to be cartesian closed.
Our models $\poly$ and $\dill$ are analogues of the categories $\mathsf{Dial}(\mathbb C)$ and $\mathsf{Dial}(\mathbb C)_!$.
The difference here is that we show that both models admit $\Pi$-types and thus in particular are cartesian closed categories.
Our introduction to the subject of Dialectica categories continues in \S\ref{sec:dialectica-category}, guided by our goal to generalize them with dependent types.
We take \emph{display map categories} as our notion of model of type theory, and review the basics in \S\ref{sec:displ-map-categ}.
This keeps us close to \cite{DEPAIVA_thesis1991} and the concept of cartesian closed category, since the latter is a category with finite products whose product projections form a class of display maps with $\Pi$-types.

In \S\ref{sec:gluing-construction} we introduce our abstract framework of \emph{fibred display map categories} (see also \cite{UEMURA_fibredfibrationcategories2017}), which are families of display map categories indexed by some other category, and also the fundamental construction underlying our main theorems, the \emph{gluing construction}.
This turns a fibred display map category over a display map category into a new display map category, by giving a class of display maps in the total family of the underlying fibration.
Where this work most differs from related work (e.g.\ \cite{SHULMAN2014}) is in the description of the type constructors: in $\poly$ and $\dill$ the $\Pi$-types and universes are not built out of their fibred versions.
That is to say, each \emph{fibre} of our fibred display map category need not support the standard version of the type constructor we wish to build.
Instead, we introduce the notions of \emph{quasifibred $\Pi$-type} and \emph{quasifibred universe} and give conditions for these to suffice to build $\Pi$-types and universes in the glued model.

Our main results are presented in \S\ref{sec:polynomial-models} and \S\ref{sec:diller-nahm-model}.
The result for $\poly$ relies on some `extensivity' assumptions, i.e.\ a particularly well-behaved kind of finite sum type.
Thus in \S\ref{sec:finite-sums} we propose and develop the notions of \emph{semi-extensive}, \emph{extensive}, and \emph{strongly extensive} display map category, the last of which contains what we need to build the dependent products.
In \S\ref{sec:biproducts} we recall some facts about additive monads and biproducts that are necessary to construct $\dill$.
These are used in \S\ref{sec:diller-nahm-model} for modelling the finite multiset types required for the Diller-Nahm variant.
While this account is self-contained, we refer the reader to \cite{VONGLEHN2014,MOSS2017} for some proofs.

\subsection*{Related work}

The `simply-typed' Dialectica categories have been considered in \cite{DEPAIVA_thesis1991,HYLAND2002,BIERING_ccdc2008}.
Topos- and tripos-theoretic versions have also been studied \cite{BBLBCB2007}.
The general gluing construction is related to the oplax limits of \cite{SHULMAN2014}, though we deal with a slightly different situation.
Essentially the same basic situation of \emph{fibred display map categories} has been considered independently by Uemura in \cite{UEMURA_fibredfibrationcategories2017} as \emph{fibred type-theoretic fibration categories}.
However, our work differs from \cite{SHULMAN2014,UEMURA_fibredfibrationcategories2017} in that our construction of $\Pi$-types and universes are more general.
The problem of finding a factorization system in the total category of a fibration has also arisen in the study of model categories, e.g.\ \cite{Stanculescu2012}.
Our gluing construction is closely related to the idea of \emph{logical predicates} or \emph{logical relations} as used for parametricity \cite{MAREYNOLDS1992,HERMIDA1993}.

\section{The Dialectica category}\label{sec:dialectica-category}

The idea of a category-theoretic version of the Dialectica interpretation was put forward by de Paiva \cite{DEPAIVA_thesis1991}.
This `Dialectica construction' takes as input some basic category $\mathbb C$, assumed at least to have finite products, and produces a new category $\mathsf{Dial}(\mathbb C)$.
An object of the Dialectica category $\mathsf{Dial}(\mathbb C)$ is a relation in $\mathbb C$, i.e. a pair $(U,X)$ of objects in $\mathbb C$ together with a monomorphism $\alpha : A \rightarrowtail U \times X$.
An arrow $(U,X,\alpha) \to (V,Y,\beta)$ is a realization of \eqref{eq:2}, i.e. it consists of arrows $f : {} U \to V$ and $F : {} U \times Y \to X$ in $\mathbb C$ such that $(\pi_U,F)^*(\alpha) \leq (f \times 1_Y)^*(\beta)$ as objects of the subobject preorder $\mathsf{Sub}(U \times Y)$.
The reader should refer to \cite{DEPAIVA_thesis1991} for a useful convention for diagrammatically depicting such arrows.

In \cite{DEPAIVA_thesis1991} it was shown that $\mathsf{Dial}(\mathbb C)$ admits a symmetric monoidal product and moreover, when $\mathbb C$ is cartesian closed, $\mathsf{Dial}(\mathbb C)$ is monoidal closed.
When $\mathbb C$ is an \emph{extensive category} (see Definition \ref{definition-extensive}), meaning that finite coproducts are well-behaved, there is a cartesian product on $\mathsf{Dial}(\mathbb C)$ given by
\begin{displaymath}
  (U,X,\alpha) \times (V,Y,\beta)
  =
  (U \times V, X + Y, \theta),
\end{displaymath}
where to form $\theta \in \mathsf{Sub}(U \times V \times (X + Y))$ we use the identification $\mathsf{Sub}(U \times V \times (X + Y)) \cong \mathsf{Sub}(U \times V \times X) \times \mathsf{Sub}(U \times V \times Y)$ which is a consequence of extensivity and then for the two components we give the two weakenings $\pi_{U,X}^*(\alpha) \in \mathsf{Sub}(U \times V \times X)$ and $\pi_{V,Y}^*(\beta) \in \mathsf{Sub}(U \times V \times Y)$ of $\alpha$ and $\beta$ along the respective product projections $\pi_{U,X} : U \times V \times X \to U \times X$ and $\pi_{V\!,Y} : U \times V \times Y \to V \times Y$.

\subsection*{The Dialectica category as a fibred category}

In \cite{HYLAND2002}, we see the beginning of a fibred approach to a more general Dialectica construction.
In its most abstract form, given a composable pair of fibrations
\vspace{-0.3cm}
\begin{equation}\label{eq:5}
  \mathbb E \xrightarrow{g} \mathbb D \xrightarrow{f} \mathbb C,
\end{equation}
the Dialectica construction amounts to taking the total category of $(f \circ (g^{\,\mathsf{op}}))^{\,\mathsf{op}}$, where `$\mathsf{op}$' denotes the \emph{opposite fibration} (or fibrewise opposite category, see \cite{STREICHER_alajeanbenabou2014}).
In \cite{HYLAND2002}, Hyland constructs $\mathsf{Dial}(p)$ by taking the $f$ above to be the simple fibration $\mathbb C_{(-)} \to \mathbb C$ of a category $\mathbb C$ with finite products.
An object over $I$ is a pair $(I,A)$ and an arrow $(I,A) \to (J,B)$ is an arrow $I \times A \to B$ in $\mathbb C$ (see \cite{JACOBS_categoricallogic1999} for details).
One supposes also a preordered fibration $p : \mathbb P \to \mathbb C$ of `predicates', which could be the subobject fibration.
Then one takes $g$ to be the pullback of $p$ along the `comprehension' $\mathbb C_{(-)} \to \mathbb C$ given by $(I,A) \mapsto I \times A$.
Applying the abstract Dialectica construction, we get a category $\mathsf{Dial}(p)$ whose objects are triples $(U,X,\alpha)$ this time with $\alpha \in \mathbb P(U \times X)$.
Moreover, it comes naturally fibred over $\mathbb C$ via the projection $(U,X,\alpha) \mapsto U$.

\subsection*{A cartesian closed Dialectica category}

While $\mathsf{Dial}(\mathbb C)$ is symmetric monoidal closed, it is not in general cartesian closed.
Our Theorem \ref{theorem-polynomial} gives a cartesian closed Dialectica category, after passage to dependent types.
The essential point is that we can take a more general fibration $f : \mathbb D \to \mathbb C$ as in \eqref{eq:5} rather than the simple fibration, which corresponds to a trivial type dependency in which the only type families are constant ones.
In terms of the Dialectica interpretation this generalization says that, in propositions of the form $\exists\, u^U \forall x^X A(u;x)$, the type $X$ of the variable $x$ is permitted to depend on the first variable, i.e.\ $X$ is a dependent type over $U$.

Let us outline why this extra generality matters, by showing it arises in the attempt to construct a function space for two objects of $\mathsf{Dial}(\mathbb C)$.
Given objects $(U,X,\alpha)$, $(V,Y,\beta)$ and $(W,Z,\gamma)$ of $\mathsf{Dial}(\mathbb C)$, the function space $(I,A,\omega) = (V,Y,\beta) \Rightarrow (W,Z,\gamma)$ needs to classify maps $(U,X,\alpha) \times (V,Y,\beta) \to (W,Z,\gamma)$ by maps $(U,X,\alpha) \to (I,A,\omega)$.
Ignoring the `subobject' part for now and focussing on the types (this is in fact taking $p : \mathbb P \to \mathbb C$ to be the terminal/identity fibration over $\mathbb C$ in the setting of \cite{HYLAND2002}), we need to classify pairs of maps $f : U \times V \to W$ and $F : U \times V \times Z \to X + Y$.
We break up $F$ as
\vspace{-0.1cm}
\begin{align*}
  F_1 :{}& U \times V \times Z \to Y + 1 \\
  F_2 :{}& (U \times V \times Z) \backslash \mathop{\mathsf{dom}} F_1 \to X
\end{align*}
where we consider $F_1$ as a partial function $U \times V \times Z \rightharpoonup Y$.
Thus we take $I$ to be the object $(V \Rightarrow W) \times (V \times Z \Rightarrow X + 1)$ and observe that we are now stuck without being able to take the object $A$ to be the dependent type
\begin{displaymath}
  \langle g,G \rangle : I \vdash (V \times Z)\backslash \mathop{\mathsf{dom}} G,
\end{displaymath}
thinking of $G : V \times Z \Rightarrow X + 1$ as a partial function.
In fact, it is not entirely straightforward to give a satisfactory meaning to this proposed object $A$ in a dependent type theory.
In order to show it has the desired property, we need a notion of display map category with well-behaved finite sum types and we will see that, in particular, the coproduct inclusions must be display maps.
We spend \S\ref{sec:finite-sums} building the theory of \emph{strongly extensive display map categories}, which we apply in \S\ref{sec:polynomial-models} to build the Dialectica model of type theory.

\section{Display map categories}\label{sec:displ-map-categ}

We will take the following as our basic notion of model of type theory, which we recall from \cite{TAYLOR1999}.
\begin{definition}\label{definition-dmc}
  A \emph{display map category} is a category $\mathbb B$ together with a class $\mathcal D$ of arrows containing the identities, called \emph{display maps}, such that given a display map $f : X \to A$ and an arbitrary map $s : B \to A$ there exists a pullback square
  \begin{equation}\label{eq:8}
    \begin{tikzpicture}[baseline={0.35*\arrow cm}]
      \node (b) at (0,0) {$B$};
      \node (a) at ({0.8*\arrow},0) {$A$};
      \node (y) at (0,{0.7*\arrow}) {$Y$};
      \node (x) at ({0.8*\arrow},{0.7*\arrow}) {$X$};

      \path[->]
      (b) edge node[auto,swap] {$s$} (a)
      (x) edge node[auto] {$f$} (a)
      (y) edge node[auto,swap] {$g$} (b)
      (y) edge node[auto] {$t$} (x)
      ;
    \end{tikzpicture}
  \end{equation}
  in $\mathbb B$ and, moreover, in any such square the map $g : Y \to B$ is also a display map. We denote display maps in diagrams with double-headed arrows $\twoheadrightarrow$.
  A \emph{class of display maps} in a category $\mathbb B$ is a class $\mathcal D$ of morphisms such that $(\mathbb B,\mathcal D)$ is a display map category.
\end{definition}
The intuition is that the category $\mathbb B$ is a category of contexts for a type theory and the class of display maps abstracts the collection of context projection morphism $\Gamma.A \to \Gamma$.
Hence, display map categories themselves give only the most basic structure for modelling dependent types.
Observe that if $(\mathbb B, \mathcal D)$ is a display map category, then $\mathcal D$ may be considered as a full subcategory of the arrow category $\mathbb B^\to$ of $\mathbb B$.
The codomain functor $\mathsf{cod} : \mathcal D \to \mathbb B$ is, by virtue of the pullback stability property of $\mathcal D$, a fibration.
As categories over $\mathbb B$, there is a cartesian functor $i_\mathcal D : \mathcal D \to \mathbb B^\to$ given by the inclusion.
Moreover, since $\mathcal D$ contains the identities (and hence isomorphisms), as a fibred category it has fibred terminals which are preserved by the inclusion $\mathcal D \to \mathbb B^\to$.
We will usually consider the following additional properties.

\begin{definition}
  A display map category $(\mathbb B,\mathcal D)$ is \emph{well-rooted} if $\mathbb B$ has a terminal object $1$ and for each object $X \in \mathbb B$, the unique arrow $X \to 1$ is in $\mathcal D$.
\end{definition}

\begin{definition}
  A display map category $(\mathbb B,\mathcal D)$ \emph{has $\Sigma$-types} (or \emph{dependent sums}) if $\mathcal D$ is closed under composition.
\end{definition}

\subsection*{Slices of display map categories}

The category $\mathcal D/I$ defined below represents the category of types in context $I$.
When $\mathcal D$ has $\Sigma$-types, this category is itself a display map category.

\begin{definition}
  Let $(\mathbb B,\mathcal D)$ be a display map category and let $I$ be an object of\/ $\mathbb B$.
  Then the \emph{display slice category over $I$} is the full subcategory $\mathcal D/I$ of the slice $\mathbb B/I$ whose objects are members of $\mathcal D$.
  We define the class $\mathcal D_I$ of morphisms in $\mathcal D/I$ to be given by those triangles whose underlying arrow is in $\mathcal D$.
\end{definition}

\begin{lemma}
  Let $(\mathbb B,\mathcal D)$ be a display map category with $\Sigma$-types and $I$ an object of\/ $\mathbb B$.
  Then $\mathcal D_I$ is a class of display maps in $\mathcal D/I$, and $(\mathcal D/I,\mathcal D_I)$ is well-rooted and has $\Sigma$-types.
\end{lemma}

\begin{lemma}
  Let $(\mathbb B,\mathcal D)$ be a display map category.
  For every map $h : I \to J$ in $\mathbb B$, the functor $h^* : \mathcal D/J \to \mathcal D/I$ given by pullback along $h$ preserves display maps, i.e.\ maps the class $\mathcal D_J$ into $\mathcal D_I$.
\end{lemma}

\subsection*{$\Pi$-types}

It will be useful for us to give a more general definition of dependent product.
We make use of the notion of \emph{couniversal arrow} \cite{MACLANE_cwm1998}, which allows us to consider a `right adjoint' to a functor which is defined only at a restricted class of objects.

\begin{definition}
  Let $(\mathbb B,\mathcal D)$ be a display map category and let $\mathcal E, \mathcal F \subseteq \mathcal D$ be two subclasses of $\mathcal D$.
  A display map category \emph{has $\mathcal E$-products of\/ $\mathcal F$-maps} if for every $f : B \twoheadrightarrow A$ in $\mathcal E$ and $g : C \twoheadrightarrow B$ in $\mathcal F$ there exists an $f^*$-couniversal arrow with codomain $g$ and domain in $\mathcal D$ where $f^* : \mathbb B/A \to \mathbb B/B$ is the pullback functor between slice categories.
\end{definition}

\begin{definition}
  A display map category $(\mathbb B,\mathcal D)$ \emph{has $\Pi$-types} (or \emph{dependent products}) if it has $\mathcal D$-products of $\mathcal D$-maps.
\end{definition}

The usual definition of $\Pi$-types is phrased in terms of $\mathsf{cod} : \mathcal D \to \mathbb B$ having \emph{$\mathcal D$-products}, meaning that for each $f : B \twoheadrightarrow A$ in $\mathcal D$, there is a right adjoint to the pullback functor $f^* : \mathcal D/A \to \mathcal D/B$, and this family of right adjoints satisfies the \emph{Beck-Chevalley condition} \cite{JACOBS_categoricallogic1999}.
It is an easy exercise to see that our definition, which does not mention a Beck-Chevalley condition, is equivalent, using the fact that $\mathcal D$ is stable under pullback.
We note that if $(\mathbb B,\mathcal D)$ has $\Sigma$-types, then $\mathsf{cod}$ has the dual structure, $\mathcal D$-sums.

\subsection*{Identity types}

Following the result of \cite{GAMBINOGARNER_identitywfs2008} that identity types give rise to a factorization system on the syntactic category, it is now standard to define identity types in terms of a factorization system, as in \cite{SHULMAN2014}.
That formulation only applies to well-rooted display map categories.
For general display map categories, it seems natural to phrase it in terms of slices.

\begin{definition}
  Let $\mathbb C$ be a category and $\mathcal F$ any class of maps.
  Then the \emph{left class of $\mathcal F$} is the class ${}^\square \mathcal F$ of \emph{left maps}, i.e.\ those maps $m : A \to B$ such that for any $f : X \to Y$ in $\mathcal F$ and any maps $u : A \to X$ and $v : B \to Y$ making the square
  \begin{displaymath}
    \begin{tikzpicture}
      \node (a) at (0,{0.75*\arrow}) {$A$};
      \node (b) at (0,0) {$B$};
      \node (x) at ({0.8*\arrow},{0.75*\arrow}) {$X$};
      \node (y) at ({0.8*\arrow},0) {$Y$};

      \path[->]
      (a) edge node[auto,swap] {$m$} (b)
      (a) edge node[auto] {$u$} (x)
      (b) edge node[auto,swap] {$v$} (y)
      (b) edge[dotted] node[auto,outer sep=-2] {$h$} (x)
      (x) edge node[auto] {$f$} (y)
      ;
    \end{tikzpicture}
  \end{displaymath}
  commute, there exists a dotted map $h$ making both triangles in that diagram commute.
\end{definition}

\begin{definition}
  A display map category $(\mathbb B,\mathcal D)$ has \emph{stable left classes} if, for any morphism $h : I \to J$ in $\mathbb B$, the functor $h^* : \mathcal D/J \to \mathcal D/I$ given by pullback along $h$ preserves left maps, i.e.\ maps the class ${}^\square\mathcal D_J$ into ${}^\square\mathcal D_I$.
\end{definition}

\begin{definition}
  Let $\mathbb C$ be a category and $\mathcal A$ and $\mathcal B$ any two classes of morphisms in $\mathbb C$.
  Then $\mathbb C$ \emph{admits\/ $(\mathcal A,\mathcal B)$-factorizations} if for any map $h : X \to Y$ in $\mathbb C$, there exists a factorization $h = fi$
  where $i \in \mathcal A$ and $f \in \mathcal B$.
\end{definition}

\begin{definition}
  A display map category $(\mathbb B,\mathcal D)$ \emph{has identity types} if it has stable left classes and for every object $I \in \mathbb B$, the category $\mathcal D/I$ admits $({}^\square\mathcal D_I,\mathcal D_I)$-factorizations.
\end{definition}

In good situations, this definition is equivalent to a `global' one.

\begin{lemma}\label{lemma-identity-types-well-rooted}
  Let $(\mathbb B,\mathcal D)$ be a well-rooted display map category with $\Sigma$-types.
  Then $(\mathbb B,\mathcal D)$ has identity types if and only if\/ $\mathbb B$ admits $({}^\square \mathcal D,\mathcal D)$-factorizations and ${}^\square\mathcal D$ is stable in the sense that for any display maps $f : A \twoheadrightarrow J$ and $g : B \twoheadrightarrow J$, any left map $m : A \to B$ for which $g \circ m = f$, and any map $h : I \to J$, the pullback of $m$ along $h$ is a left map.
  That is to say, in the diagram
  \begin{displaymath}
    \begin{tikzpicture}
      \node (i) at (0,0) {$I$};
      \node (j) at (3,-0.4) {$J$};
      \node (a) at (2.3,0.7) {$A$};
      \node (b) at (3.7,0.9) {$B$};
      \node (ap) at (-0.7,1.1) {$A'$};
      \node (bp) at (0.7,1.3) {$B'$};

      \path[->]
      (a) edge node[auto,swap,pos=0.3] {$m$} (b)
      (i) edge node[auto,swap] {$h$} (j)
      (ap) edge (a)
      (bp) edge (b)
      (ap) edge[dotted] node[auto] {$\overline m$} (bp)
      ;
      \path[->>]
       (a) edge node[auto,swap,pos=0.3] {$f$} (j)
      (b) edge node[auto,pos=0.3] {$g$} (j)
      (ap) edge node[auto,swap,pos=0.3] {$f'$} (i)
      (bp) edge node[auto,pos=0.3] {$g'$} (i)
      ;
    \end{tikzpicture}
  \end{displaymath}
  where both near-vertical squares are pullbacks, the unique dotted arrow $\overline m$ making the diagram commute is a left map.
\end{lemma}

In fact, we can do slightly better.
If $(\mathbb B,\mathcal D)$ is well-rooted and has $\Sigma$- and $\Pi$-types, then the stability of ${}^\square \mathcal D$ in the sense of Lemma \ref{lemma-identity-types-well-rooted} follows from the existence of $({}^\square \mathcal D,\mathcal D)$-factorizations \cite[Lemma 7.2]{SHULMAN_eidiagrams2017}.
We will not make use of this simplification since we wish to consider the various type constructors separately.
Another point to make is that our definition of identity types is really slightly too weak.
We should include a condition that left maps are stable under pullback along display maps.
This issue is discussed in \cite[\S3.4.3]{LUMSDAINEWARREN2014}, for example, but it disappears in the presence of $\Pi$-types.

\begin{lemma}\label{lemma-identity-and-pi}
  Let $(\mathbb B,\mathcal D)$ be a display map category with $\Pi$-types.
  Then for each $I$,\/ ${}^\square (\mathcal D_I)$ is stable under pullback along $\mathcal D_I$ in $\mathcal D/I$.
\end{lemma}

\subsection*{Universes}

\begin{definition}
  A \emph{universe} in a display map category $(\mathbb B,\mathcal D)$ is a display map $u : \widetilde{\mathcal U} \to \mathcal U$ such that the class $\mathcal D_u$ of maps which are pullbacks of $u$ is a (not necessarily well-rooted) class of display maps in $\mathbb B$.
  Such maps are called \emph{$u$-small}.
  The universe \emph{has $\Sigma$-types} if $(\mathbb B,\mathcal D_u)$ has $\Sigma$-types.
  The universe \emph{has $\Pi$-types} if $(\mathbb B,\mathcal D_u)$ has $\Pi$-types.
  Equivalently, this says that $(\mathbb B,\mathcal D)$ has $\mathcal D_u$-products of $\mathcal D_u$-maps which are again in $\mathcal D_u$.
  In the situation where $(\mathbb B,\mathcal D)$ has identity types, we say that the universe \emph{has identity types} if for every object $I \in \mathbb B$, the category $\mathcal D_u/I$ admits $({}^\square \mathcal D_I,(\mathcal D_u)_I)$-factorizations.
  (Note that since ${}^\square ((\mathcal D_u)_I)$ may be larger than ${}^\square(\mathcal D_I)$, and need not be stable even if the latter is, this is not the same as asking for $(\mathbb B,\mathcal D_u)$ to have identity types).
  When $(\mathbb B,\mathcal D)$ has finite sum types (see Definition \ref{definition-finite-sum-types}), the universe is \emph{closed under finite sums} if $(\mathbb B,\mathcal D_u)$ has finite sum types which are preserved by the inclusion into $(\mathbb B,\mathcal D)$.
\end{definition}

\section{The gluing construction}\label{sec:gluing-construction}

We give details of the basic general construction used to build the models of \S\ref{sec:polynomial-models} and \S\ref{sec:diller-nahm-model}, which we call \emph{gluing}.
This is closely related to the work of Shulman \cite{SHULMAN2014} on oplax limits of models of type theory: our gluing construction is a different generalization of the special case referred to as the `gluing construction' there.

In our setting we start with some base model of type theory given by a display map category $(\mathbb B,\mathcal D)$.
Then we have a system of `new types' for each context in $\mathbb B$, with which we want to extend our original type theory.
This is modelled by a fibration $p : \mathbb E \to \mathbb B$ for which each fibre category $\mathbb E(I)$ is itself a display map category --- we call this a \emph{fibred display map category}.
The result of the gluing construction is a class of display maps in $\mathbb E$ making $\mathbb E$ into a model of type theory.
We investigate which type constructors exist in $\mathbb E$ given certain hypotheses on $(\mathbb B,\mathcal D)$ and the fibred display map category $(p : \mathbb E \to \mathbb B,\mathcal E)$.

The construction of \cite{SHULMAN2014} covers the case where the fibration $p$ arises in the style of `logical relations' \cite{MAREYNOLDS1992}, i.e.\ as the pullback of the self-indexing along some morphism of models $\mathbb E \to \mathbb C$.
Our more general situation has been considered independently in \cite{UEMURA_fibredfibrationcategories2017}.

\subsection*{Fibred display map categories}

The following definition corresponds to that of `fibred type-theoretic fibration category' in \cite{UEMURA_fibredfibrationcategories2017} but with only the basic structure of dependent type theory (i.e.\ no type constructors).

\begin{definition}
  A \emph{fibred display map category} consists of a fibration $p : \mathbb E \to \mathbb B$ together with, for each object $I \in \mathbb B$, a class $\mathcal E_I$ of display maps in the fibre category $\mathbb E(I)$ such that, for every arrow $h : I \to J$ in $\mathbb B$, the action of reindexing $h^* : \mathbb E(J) \to \mathbb E(I)$ preserves display maps, i.e.\ it maps the class $\mathcal E_J$ into $\mathcal E_I$.
  Moreover, each $h^*$ must preserve all pullbacks of display maps.
\end{definition}

We can also collect all of the `fibrewise' classes of display maps into one class $\mathcal E = \bigcup_{I \in \mathbb B} \mathcal E_I$ of $p$-vertical display maps.
This very nearly makes $(\mathbb E,\mathcal E)$ into a display map category except that, for example, $\mathcal E$ need not contain all of the isomorphisms.
As an aside, we note that there is a natural generalization of display map categories which encompasses it, where from Definition \ref{definition-dmc} we only require the $g$ in \eqref{eq:8} to be in $\mathcal D$ for at least one choice of pullback square rather than for every pullback square.
The theory of such structures is largely unchanged, and they are naturally seen as a special case of full comprehension categories \cite{JACOBS_compcats1993}.

\begin{definition}
  Let $p : \mathbb E \to \mathbb B$ be a fibration over a display map category $(\mathbb B,\mathcal D)$.
  Then we define the class $\overline{\mathcal D}$ in $\mathbb E$ to consist of those $p$-cartesian morphisms $f$ for which $p(f) \in \mathcal D$.
\end{definition}

We note, as an aside, that this time $(\mathbb E,\overline{\mathcal D})$ \emph{is} a display map category.
The following is straightforward.

\begin{proposition}\label{proposition-glued-dmc}
  Let $(\mathbb B,\mathcal D)$ be a display map category and let $(p : \mathbb E \to \mathbb B,\mathcal E)$ be a fibred display map category.
  Define the class $\mathcal G$ to be those morphisms $f$ in $\mathbb E$ with cartesian component in $\overline{\mathcal D}$ and vertical component in $\mathcal E$.
  Then $\mathcal G$ is a class of display maps in $\mathbb E$.
  Moreover, if\/ $(\mathbb B,\mathcal D)$ and each $(\mathbb E(I),\mathcal E_I)$ are well-rooted, then so is $(\mathbb E,\mathcal G)$.
\end{proposition}

\begin{definition}
The \emph{glued display map category} is the category $\mathbb E$ together with the class of display maps $\mathcal G$.
\end{definition}

In \cite{SHULMAN2014,UEMURA_fibredfibrationcategories2017} the members of $\mathcal G$ are referred to as \emph{Reedy fibrations}.
We continue with the notation of Proposition \ref{proposition-glued-dmc} throughout \S\ref{sec:gluing-construction}.

\subsection*{Type constructors in a glued display map category}

In general, type constructors in a glued display map category can be built out of type constructors in the base and each fibre display map category, given sufficient compatibility conditions.
The model example for such constructions is the following well-known result.

\begin{proposition}[{\cite[Corollary 4.9]{HERMIDA_2fibs1999}}]\label{proposition-hermida-finite-fibred-products}
  Let $q : \mathbb D \to \mathbb C$ be a fibration over a base with finite products.
  Then $q$ has fibred finite products if and only if\/ $\mathbb D$ has finite products preserved by $q$.
\end{proposition}

\subsection*{$\Sigma$-types}

\begin{definition}
  We say that a fibred display map category $(p : \mathbb E \to \mathbb B,\mathcal E)$ has \emph{fibrewise $\Sigma$-types} if, for each object $I \in \mathbb B$, the class $\mathcal E_I$ is closed under composition.
  Equivalently, if the class $\mathcal E$ is closed under composition.
\end{definition}

The following is straightforward \cite[Proposition 3.12]{VONGLEHN2014}.

\begin{proposition}\label{proposition-fibred-sigma}
  Suppose that $(p : \mathbb E \to \mathbb B,\mathcal E)$ has fibred $\Sigma$-types, and that $(\mathbb B,\mathcal D)$ has $\Sigma$-types.
  Then $(\mathbb E,\mathcal G)$ has $\Sigma$-types.
\end{proposition}

\subsection*{Identity types}

For the construction of identity types we assume that $\mathbb B$ and all of the fibre display map categories $(\mathbb E(I),\mathcal E_I)$ are well-rooted and have $\Sigma$-types, hence we may use the characterization of identity types from Lemma \ref{lemma-identity-types-well-rooted}.

Let us first identify the left class ${}^\square\mathcal G$, since it leads us to what seems to be a fundamental condition on a fibration over a display map category.
It is easy to check by the basic properties of left classes that ${}^\square\mathcal G = {}^\square\mathcal E \cap {}^\square\overline{\mathcal D}$.
The next lemma follows easily from properties of fibrations.

\begin{lemma}
  The class\/ ${}^\square\overline{\mathcal D}$ consists of precisely those\/ $i : A \to B$ in $\mathbb E$ such that $p(i) \in {}^\square\mathcal D$.
\end{lemma}

\begin{lemma}
  The class of vertical maps in\/ ${}^\square\mathcal E$ (and hence of those in ${}^\square\mathcal G$) is precisely\/ $\bigcup_{I \in \mathbb B} {}^\square\mathcal E_I$.
\end{lemma}

We cannot get further without an additional assumption.
The left map condition appears as condition 3) in \cite[Theorem 3.3]{UEMURA_fibredfibrationcategories2017}.

\begin{definition}
  \emph{The left map condition} says that for any $m : I \to J$ in ${}^\square \mathcal D$ and any $X \in \mathbb E(J)$, the reindexing functor $m^*/X : \mathcal E_J/X \to \mathcal E_I/m^*(X)$, between slices of fibre categories, is full.
\end{definition}

Using closure of $\mathcal E_J$ under pullbacks in $\mathbb E(J)$, to check the left map condition it is sufficient to check that for any $\mathcal E_J$-map $f : Y \twoheadrightarrow X$, the sections of $f$ in $\mathbb E(J)$ are mapped surjectively to the sections of $m^*(f)$ in $\mathbb E(I)$ by $m^*$.
A related formulation is the following.

\begin{lemma}
  The left map condition holds if and only if every $p$-cartesian map $f$ with $p(f) \in {}^\square\mathcal D$ is in ${}^\square \mathcal E$.
\end{lemma}

\begin{proposition}[{\cite[Proposition 4.6.4]{MOSS2017}}]\label{proposition-left-map-condition-left-class}
  The left map condition holds if and only if\/ ${}^\square\mathcal G$ consists of all maps lying over a map in ${}^\square\mathcal D$ with vertical component in $\bigcup_{I \in \mathbb B} {}^\square \mathcal E_I$.
\end{proposition}

\begin{definition}
  A fibred display map category $(p : \mathbb E \to \mathbb B, \mathcal E)$ has \emph{fibrewise identity types} if, for each object $I \in \mathbb B$, $(\mathbb E(I),\mathcal E_I)$ has identity types and the fibrewise left classes are stable under reindexing, i.e.\ for any map $h : I \to J$, $h^* : \mathbb E(J) \to \mathbb E(I)$ maps the class ${}^\square \mathcal E_J$ into ${}^\square \mathcal E_I$.
\end{definition}

\begin{theorem}\label{theorem-fibred-identity-types}
  Suppose that $(\mathbb B,\mathcal D)$ has identity types and that $(p : \mathbb E \to \mathbb B,\mathcal E)$ has fibred identity types and satisfies the left map condition.
  Then $(\mathbb E,\mathcal G)$ has identity types.
\end{theorem}
\begin{proof}
  The factorization of an arrow $f : B \to A$ in $\mathbb E$ is the construction given in \cite[Proposition 3.21]{VONGLEHN2014} (but see also \cite[Lemma 3.9]{UEMURA_fibredfibrationcategories2017}).
  It may be read from Figure \ref{fig:identity-factorization} as the ${}^\square\mathcal G$-map $x \circ \overline \imath : B \to K_v$ composed with the $\mathcal G$-map $\overline g \circ y : K_v \to A$.
  This is constructed by factorizing $p(f)$ in $\mathbb B$ as $g \circ i$, and using well-rootedness of $\mathcal D$ to construct a retraction $r$ of $i$.
  Then $\overline r$ is taken to be a cartesian lift of $r$ with codomain $B$, whence $\overline \imath$ is a cartesian lift of $i$ with codomain $Q$, which may be chosen to have domain $B$ since $ri = 1_{p(B)}$.
  Taking $\overline g$ to be a cartesian lift of $g$ with codomain $A$, we construct $l : B \to P$ as the factorization of $f$ through $\overline g$ lying over $i$.
  Now $v : Q \to P$ is taken to be a filler for the square
  \begin{displaymath}
    \begin{tikzpicture}
      \node (q) at (0,0) {$Q$};
      \node (b) at (0,{0.75*\arrow}) {$B$};
      \node (a) at ({0.8*\arrow},0) {$\top_{K_{p(f)}}$};
      \node (p) at ({0.8*\arrow},{0.75*\arrow}) {$P$};

      \path[->]
      (b) edge node[auto,swap] {$\overline\imath$} (q)
      (q) edge (a)
      (b) edge node[auto] {$l$} (p)
      (p) edge (a)
      (q) edge[dotted] node[auto,outer sep=-0.05cm] {$v$} (p)
      ;
    \end{tikzpicture}
  \end{displaymath}
  where $\top_{K_{p(f)}}$ is the terminal object in $\mathbb E(K_{p(f)})$.
  It follows that $v$ is a vertical map, and hence we can use the identity types in $\mathbb E(K_{p(f)})$ to factorize it as $y \circ x$.

  By Proposition \ref{proposition-left-map-condition-left-class}, $x \circ \overline\imath \in {}^\square \mathcal G$, and clearly $\overline g \circ y \in \mathcal G$.
  Only the stability condition of \ref{lemma-identity-types-well-rooted} remains.
  This is easy to verify, and a proof can be found in \cite[Lemma 6.4.8]{MOSS2017}.
  Alternatively, as we remarked before, stability follows if we have $\Sigma$- and $\Pi$-types.
\end{proof}

\begin{figure}
  \centering
  \begin{displaymath}
    \begin{tikzpicture}
      \node (b) at (0,0) {$B$};
      \node (a) at (4,0) {$A$};

      \node (q) at (1.3,0.5) {$Q$};
      \node (p) at (2.7,0.5) {$P$};
      \node (kv) at (2.05,0.9) {$K_v$};

      \node (pb) at (0,-2) {$p(B)$};
      \node (pa) at (4,-2) {$p(A)$};
      \node (kpf) at (2,-1.3) {$K_{p(f)}$};

      \node (ee) at (5.5,0.5) {$\mathbb E$};
      \node (bb) at (5.5,-1.8) {$\mathbb B$};

      \path[->]
      (ee) edge node[auto] {$p$} (bb)
      ;

      \path[white!50!black,dotted]
      (b) edge (pb)
      (a) edge (pa)
      (q) edge (kpf)
      (p) edge (kpf)
      (kv) edge (kpf)
      ;

      \path[->]
      (pb) edge[out=-15,in=-165] node[auto,swap] {$p(f)$} (pa)
      (b) edge[out=-15,in=-165] node[auto,swap] {$f$} (a)
      (pb) edge node[auto,swap,outer sep=-1] {$i$} (kpf)
      (kpf) edge[->>] node[auto] {$g$} (pa)
      (kpf) edge[out=-180,in=40] node[auto,swap] {$r$} (pb)
      ;

      \path[->]
      (q) edge[out=-180,in=40] node[auto,swap] {$\overline r$} (b)
      (b) edge node[auto,swap,outer sep=-2,pos=0.7] {$\overline \imath$} (q)
      (p) edge[->>] node[auto] {$\overline g$} (a)
      (q) edge node[auto,swap,pos=0.4] {$v$} (p)
      (q) edge[auto,shorten >={-0.08cm}] node {$x$} (kv)
      (kv) edge[->>,shorten <={-0.1cm}] node[auto] {$y$} (p)
      ;

      \path[->]
      (b) edge[out=-5,in=-150] node[auto,swap,pos=0.85,outer sep=-1.5] {$l$} (p);
      ;
    \end{tikzpicture}
  \end{displaymath}
  \caption{The identity factorization of $f : B \to A$.}
  \label{fig:identity-factorization}
\end{figure}

\subsection*{$\Pi$-types}

There is a result similar to Proposition \ref{proposition-hermida-finite-fibred-products} showing the equivalence of fibrewise cartesian closed structure with ordinary cartesian closed structure in the total category when the fibration has simple products \cite[Corollary 4.12]{HERMIDA_2fibs1999}.
This generalizes from the simply-typed case to the dependently-typed case.
The following result connecting $\Pi$-types in the glued model to fibrewise $\Pi$-types is \cite[Proposition 3.14]{VONGLEHN2014} (but see also \cite[Lemma 3.10]{UEMURA_fibredfibrationcategories2017}).

\begin{proposition}\label{proposition-fibred-pi}
  Suppose that $(\mathbb B, \mathcal D)$ has $\Pi$-types.
  Then $(\mathbb E, \mathcal G)$ has $\Pi$-types preserved by $p : \mathbb E \to \mathbb B$ if and only if each fibre category $(\mathbb E(I), \mathcal E_I)$ has $\Pi$-types which are stable under reindexing and the fibration $p$ has $\mathcal D$-products which preserve the display maps in $\mathcal E$.
\end{proposition}

However, in our Dialectica models, $\Pi$-types in $(\mathbb E, \mathcal G)$ are not preserved by the fibration $p$.
But the fibred display map category will have $\mathcal D$-products which preserve display maps.
With the following, we can consider $\overline{\mathcal D}$- and $\mathcal E$-products separately.

\begin{lemma}\label{lemma-pi-type-d-products}
  Suppose that $(\mathbb B,\mathcal D)$ has $\Pi$-types.
  Then $p : \mathbb E \to \mathbb B$ has $\mathcal D$-products which preserve the fibrewise display maps if and only if\/ $(\mathbb E,\mathcal G)$ has $\overline{\mathcal D}$-products of $\mathcal G$-maps which are sent to $\mathcal D$-products of $\mathcal D$-maps by $p$.
\end{lemma}

\begin{proof}
  This is a straightforward matter of comparing the definitions, and is proved as stated in \cite[Proposition 6.5.5]{MOSS2017}.
  See also \cite[Lemma 3.13]{VONGLEHN2014}, \cite[Theorem 8.8]{SHULMAN2014}, and \cite[Lemma 3.10]{UEMURA_fibredfibrationcategories2017}.
\end{proof}

\begin{lemma}\label{lemma-e-products-of-g-maps}
  Suppose that $(\mathbb B,\mathcal D)$ has $\Sigma$-types.
  Then $(\mathbb E,\mathcal G)$ has $\mathcal E$-products of $\mathcal G$-maps if and only if it has $\mathcal E$-products of $\mathcal E$-maps.
\end{lemma}
\begin{proof}[Proof (sketch)]
  Since $\mathcal E \subseteq \mathcal G$, the only if direction is trivial.
  For the if direction, suppose that we have a $\mathcal G$-map given as a composite $vw$ where $(w : W \twoheadrightarrow Z) \in \mathcal E$ and $(v : Z \twoheadrightarrow Y) \in \overline{\mathcal D}$, and also a map $\phi : Y \to X$ in $\mathcal E$.
  The key point is that we can take a cartesian lift of $p(v)$ with codomain $X$ to get $v' : Z' \twoheadrightarrow X$ and get an induced $\mathcal E$-map $\phi' : Z \twoheadrightarrow Z'$ satisfying $v' \circ \phi' = \phi \circ v$, since $\mathcal E$ is stable under reindexing.
  Then the product of $v \circ w$ along $\phi$ is given by $v'$ composed with the product of $w$ along $\phi'$, which exists and is a $\mathcal G$-map by hypothesis.
  The condition of $\Sigma$-types in $\mathcal D$ ensures that the composition of a $\mathcal G$-map with a $\overline{\mathcal D}$-map is again in $\mathcal G$.
  Full details can be found in \cite[Proposition 6.5.7]{MOSS2017}.
\end{proof}

To connect this to Proposition \ref{proposition-fibred-pi}, observe that fibrewise $\Pi$-types exist if and only if $(\mathbb E,\mathcal G)$ has $\mathcal E$-products of $\mathcal E$-maps which are again in $\mathcal E$ (see \cite[Proposition 6.5.8]{MOSS2017}).
However, in our situation we do not have $\Pi$-types in the fibres.
Instead, we use the following.

\begin{definition}\label{definition-quasifibred-pi}
A fibred display map category $(p : \mathbb E \to \mathbb B,\mathcal E)$ with class $\mathcal D$ of display maps in the base has \emph{quasifibred\/ $\Pi$-types} if for any $I \in \mathbb B$ and composable pair of $\mathcal E_I$-maps $\psi : Z \twoheadrightarrow Y$ and $\phi : Y \twoheadrightarrow X$ in $\mathbb{E}(I)$ there exists a $\mathcal D$-map $q : Q \twoheadrightarrow I$ in $\mathbb B$ and an $\mathcal E_Q$-map $\pi : P \twoheadrightarrow q^*(X)$ in $\mathbb{E}(Q)$ together with, for every map $w : W \to X$ in $\mathbb{E}(I)$ a bijection, natural in $W$, between the set of maps $\phi^*(w) \to \psi$ in the slice $\mathbb{E}(I)/Z$ and the set of pairs $(s, t)$ where $s : I \to Q$ is a section of $q$ and $t$ is a map $w \to s^*(\pi)$ in the slice $\mathbb{E}(I)/Y$. Moreover this data must be stable under reindexing.
\end{definition}

The following lemma is straightforward after unfolding definitions, whence the following theorem is immediate.

\begin{lemma}
  The fibred display map category $(p,\mathcal E)$ over $(\mathbb B,\mathcal B)$ has quasifibred $\Pi$-types if and only if it has $\mathcal E$-products of $\mathcal E$-maps.
\end{lemma}

\begin{theorem}\label{theorem-quasifibred-pi}
  Let $(p : \mathbb E \to \mathbb B, \mathcal E)$ be a fibred display map category over a display map category $(\mathbb B,\mathcal D)$ with $\Sigma$-types.
  Then it has quasifibred dependent products and $p$ has $\mathcal D$-products which preserve the fibrewise display maps if and only if $(\mathbb E,\mathcal G)$ has $\Pi$-types such that products along $\overline{\mathcal D}$-maps are preserved by $p$.
\end{theorem}

\subsection*{Universes}

\begin{definition}
  A \emph{quasifibred universe} in a fibred display map category $(p : \mathbb E \to \mathbb B, \mathcal E)$ consists of an object $\Omega \in \mathbb B$ together with an $\mathcal E_\Omega$-map $v : \widetilde{\mathcal V} \twoheadrightarrow \mathcal V$.
  Then for each $I \in \mathbb B$, an $\mathcal E_I$-map $\phi : Y \twoheadrightarrow X$ is \emph{$v$-small} if there exists a morphism $f : I \to \Omega$ in $\mathbb B$ such that $\phi$ arises as a pullback of $f^*(v)$.
  We require of any quasifibred universe $(\Omega,v)$ that the $p$-vertical isomorphisms of $\mathbb E$ be $v$-small.
\end{definition}

The following corresponds to \cite[Proposition 4.3]{UEMURA_fibredfibrationcategories2017}, but there $\Omega$ was required to be a terminal object.

\begin{lemma}\label{lemma-quasifibred-universe}
  Let $(p : \mathbb E \to \mathbb B,\mathcal E)$ be a fibred display map category over a display map category $(\mathbb B,\mathcal D)$ with $\Pi$-types and suppose that $p$ has $\mathcal D$-products.
  Let $(\Omega,v : \widetilde{\mathcal V} \twoheadrightarrow \mathcal V)$ be a quasifibred universe in $(p,\mathcal E)$ and $u : \widetilde{\mathcal U} \twoheadrightarrow \mathcal U$ a universe in $(\mathbb B,\mathcal D)$.
  Then there exists a universe $w : \widetilde{\mathcal W} \twoheadrightarrow \mathcal W$ in $(\mathbb E,\mathcal G)$ for which the $w$-small maps are precisely those $\mathcal G$-maps which lie over a $u$-small map and whose vertical component is $v$-small.
\end{lemma}
\begin{proof}[Proof (sketch)]
  Writing $\overline u$ for the left hand map in the pullback square in $\mathbb B$
  \vspace{-0.3cm}
  \begin{displaymath}
    \begin{tikzpicture}
      \node (a) at (0,{0.8*\arrow}) {$\widetilde{\mathcal U} \times_{\mathcal U} (\widetilde{\mathcal U} \Rightarrow_{\mathcal U} \mathcal U^*(\Omega))$};
      \node (b) at ({2*\arrow},{0.8*\arrow}) {$\widetilde{\mathcal U}$};
      \node (c) at (0,0) {$\widetilde{\mathcal U} \Rightarrow_{\mathcal U} \mathcal U^*(\Omega)$};
      \node (d) at ({2*\arrow},0) {$\mathcal U$};

      \path[->>]
      (a) edge (b)
      (c) edge (d)
      (a) edge node[auto,swap] {$\overline u$} (c)
      (b) edge node[auto] {$u$} (d)
      ;
    \end{tikzpicture}
  \end{displaymath}
  where the bottom row is the exponential in $\mathcal D/\mathcal U$ of $u$ into $\pi_{\mathcal U} : \mathcal U \times \Omega \to \mathcal U$, and writing $\mathsf{ev} : \widetilde{\mathcal U} \times (\widetilde{\mathcal U} \Rightarrow_{\mathcal U} \mathcal U^*(\Omega)) \to \mathcal U \times \Omega$ for the counit of the fibred exponential, we let $\mathcal W = \Pi_{\overline u}((\pi_\Omega \circ \mathsf{ev})^*(\mathcal V))$, which lies over $\widetilde{\mathcal U} \Rightarrow_{\mathcal U} \mathcal U^*(\Omega)$.
  Now $w : \widetilde{\mathcal W} \twoheadrightarrow \mathcal W$ should lie over $\overline u$ and have vertical component the pullback of $v : \widetilde{\mathcal V} \twoheadrightarrow \mathcal V$ along
  \begin{displaymath}
    \overline u^* \Pi_{\overline u} ((\pi_\Omega \circ \mathsf{ev})^*(\mathcal V)) \to (\pi_\Omega \circ \mathsf{ev})^*(\mathcal V) \to \mathcal V
  \end{displaymath}
  where the first arrow is the counit of $\overline u \vdash \Pi_{\overline u}$.
  Further details are in \cite[Lemma 6.6.2]{MOSS2017}.
\end{proof}

The following is now easy to check from the constructions of the type constructors we have given above.
\begin{theorem}\label{theorem-quasifibred-universe}
  Let $(\mathbb B,\mathcal D)$ be a well-rooted display map category with $\Sigma$-, $\Pi$-, and identity types, and let $(p : \mathbb E \to \mathbb B, \mathcal E)$ be a fibred display map category with well-rooted fibres and fibrewise $\Sigma$-types, and satisfying the conditions of Theorems \ref{theorem-fibred-identity-types} and \ref{theorem-quasifibred-pi}.
  Suppose moreover that $(\mathbb B,\mathcal D)$ admits a universe $u : \widetilde{\mathcal U} \twoheadrightarrow \mathcal U$ closed under $\Sigma$-, $\Pi$-, and identity types, and that $(p,\mathcal E)$ admits a quasifibred universe $(\Omega,v : \widetilde{\mathcal V} \twoheadrightarrow \mathcal V)$, and let $w : \widetilde{\mathcal W} \twoheadrightarrow \mathcal W$ be the universe constructed in \ref{lemma-quasifibred-universe}.
  Suppose that for each $I \in \mathbb E$, the $v$-small maps in $\mathbb E(I)$ are closed under $\Sigma$- and identity types, and also closed under $\mathcal D$-product along $u$-small maps, and moreover $(p,\mathcal E,\mathcal D)$ admits quasifibred $\Pi$-types in such a way that whenever the $\phi$ and $\psi$ of Definition \ref{definition-quasifibred-pi} are $v$-small then the $\pi$ is $v$-small and the $q$ is $u$-small.
  Then $w : \widetilde{\mathcal W} \twoheadrightarrow \mathcal W$ is closed under $\Sigma$-, $\Pi$-, and identity types.
\end{theorem}

\section{Finite sums}\label{sec:finite-sums}

Recall from the introduction that to construct function spaces in the Dialectica category, we were led to consider a system of dependent types with the facility for forming the type family over a type of partial functions corresponding to the complements of the domains of those partial functions.
Below we propose the notion of \emph{strongly extensive finite sums} which will serve for this purpose.
It is based on the notion of \emph{extensive category}, which we recall here. Extensivity is a standard property of `categories of sets' such as any topos, as well as many `geometric' categories such as topological spaces.

\begin{definition}[\cite{CARBONILACKWALTERS_extensivedistributive1993}]\label{definition-extensive}
  A category $\mathbb C$ with finite coproducts is \emph{extensive} if in any diagram of the form
  \begin{equation}\label{eq:10}
    \begin{tikzpicture}[baseline={0.35*\arrow cm}]
      \node (a) at (0,0) {$A$};
      \node (ab) at ({\arrow},0) {$A+B$};
      \node (b) at ({2*\arrow},0) {$B$};
      \node (x) at (0,{0.75*\arrow}) {$X$};
      \node (c) at ({\arrow},{0.75*\arrow}) {$C$};
      \node (y) at ({2*\arrow},{0.75*\arrow}) {$Y$};

      \path[->]
      (x) edge (c)
      (y) edge (c)
      (a) edge (ab)
      (b) edge (ab)
      (x) edge (a)
      (c) edge (ab)
      (y) edge (b)
      ;
    \end{tikzpicture}
  \end{equation}
  where the bottom row is a coproduct diagram, we have that the top row is a coproduct diagram if and only if both squares are pullbacks.
\end{definition}

\begin{definition}\label{definition-finite-sum-types}
  A display map category $(\mathbb B,\mathcal D)$ has \emph{finite sum types} if the fibration $\mathsf{cod} : \mathcal D \to \mathbb B$ has fibred finite coproducts.
\end{definition}

Recall that a \emph{strict initial object} in a category $\mathbb B$ is an initial object $0$ such that every map of the form $X \to 0$ is an isomorphism.

\begin{proposition}\label{proposition-weak-finite-sums}
  A well-rooted display map category $(\mathbb B,\mathcal D)$ has finite sum types if and only if\/ $\mathbb B$ has finite coproducts including a strict initial object such that the copairing preserves $\mathcal D$ and commutes with the pullback of $\mathcal D$-maps, i.e.\ such that if we are given two pullback squares
  \vspace{-0.3cm}
  \begin{displaymath}
    \begin{tikzpicture}
      \node (a) at (0,{0.8*\arrow}) {$A'$};
      \node (b) at ({\arrow},{0.8*\arrow}) {$A$};
      \node (c) at (0,0) {$D$};
      \node (d) at ({\arrow},0) {$C$};

      \path[->]
      (a) edge node[auto] {$h_A$} (b)
      (c) edge node[auto,swap] {$h$} (d)
      ;
      \path[->>]
      (a) edge node[auto,swap] {$f'$} (c)
      (b) edge node[auto] {$f$} (d)
      ;

      \node (a) at ({2*\arrow},{0.8*\arrow}) {$B'$};
      \node (b) at ({3*\arrow},{0.8*\arrow}) {$B$};
      \node (c) at ({2*\arrow},0) {$D$};
      \node (d) at ({3*\arrow},0) {$C$};

      \path[->]
      (a) edge node[auto] {$h_B$} (b)
      (c) edge node[auto,swap] {$h$} (d)
      ;
      \path[->>]
      (a) edge node[auto,swap] {$g'$} (c)
      (b) edge node[auto] {$g$} (d)
      ;
    \end{tikzpicture}
  \end{displaymath}
  where $f$ and $g$ are display maps, then $[f,g] : A + B \to C$ is also a display map and the square
  \vspace{-0.3cm}
  \begin{displaymath}
    \begin{tikzpicture}
      \node (a) at (0,{0.8*\arrow}) {$A' + B'$};
      \node (b) at ({1.8*\arrow},{0.8*\arrow}) {$A + B$};
      \node (c) at (0,0) {$D$};
      \node (d) at ({1.8*\arrow},0) {$C$};

      \path[->]
      (a) edge node[auto] {$h_A + h_B$} (b)
      (c) edge node[auto,swap] {$h$} (d)
      ;
      \path[->>]
      (a) edge node[auto,swap] {$[f',g']$} (c)
      (b) edge node[auto] {$[f,g]$} (d)
      ;
    \end{tikzpicture}
  \end{displaymath}
  is a pullback.
  In this case, the sum of $f : A \twoheadrightarrow C$ and $g : B \twoheadrightarrow C$ in $\mathcal D/C$ is $[f,g] : A + B \twoheadrightarrow C$.
\end{proposition}
\begin{proof}[Proof (sketch)]
  The key point is that in the commuting square of categories and functors
  \vspace{-0.2cm}
  \begin{displaymath}
    \begin{tikzpicture}
      \node (dx) at (0,{0.8*\arrow}) {$\mathcal D/X$};
      \node (d1) at (0,0) {$\mathcal D/1$};
      \node (bx) at ({1.5*\arrow},{0.8*\arrow}) {$\mathbb B/X$};
      \node (b1) at ({1.5*\arrow},0) {$\mathbb B/1 \cong \mathbb B$};

      \path[->]
      (dx) edge node[auto] {$i_{\mathcal D,X}$} (bx)
      (d1) edge node[auto,swap ] {$i_{\mathcal D,1}$} (b1)
      (dx) edge node[auto,swap] {$\Sigma_X$} (d1)
      (bx) edge node[auto] {$\mathsf{dom}$} (b1)
      ;
    \end{tikzpicture}
  \end{displaymath}
  the functor $\Sigma_X$ is a left adjoint and $i_{\mathcal D,1}$ is an isomorphism, hence both preserve colimits, and the functor $\mathsf{dom}$ creates colimits, hence $i_{\mathcal D,X}$ also preserves colimits that exist in $\mathcal D/X$.
\end{proof}

\subsection*{Semi-extensivity}

Semi-extensive finite sum types are those for which terms in a finite coproduct context are equivalently given by terms in each of the summand contexts.

\begin{definition}
  A display map category $(\mathbb B,\mathcal D)$ is \emph{semi-extensive} if it has finite sum types and for each cospan of display maps $f : A \twoheadrightarrow I$, $g : B \twoheadrightarrow I$, the functor
  \begin{equation} \label{eq:ext}
    \mathcal D/(\mathop{\mathsf{dom}}(f +_I g)) \to \mathcal D/A \times \mathcal D/B,
  \end{equation}
  induced by reindexing along the coproduct inclusions of the coproduct of $f$ and $g$ in $\mathcal D/I$, is full and faithful.
\end{definition}

As for finite sums, we can give an equivalent `global' definition in the case of well-rooted display map categories.
\begin{proposition}
  A well-rooted display map category $(\mathbb B,\mathcal D)$ with finite sum types is semi-extensive if and only if, in any diagram of the form \eqref{eq:10}
  where the vertical arrows are display maps and the bottom row is a coproduct diagram, if both squares are pullbacks then the top row is a coproduct diagram.
\end{proposition}
\begin{proof}[Proof (sketch)]
  We merely comment that, in light of Proposition \ref{proposition-weak-finite-sums}, the definition of semi-extensive finite sums becomes much simpler in the well-rooted case.
  It says simply that for any two objects $A$ and $B$, the functor $\mathcal D/(A+B) \to \mathcal D/A \times \mathcal D/B$ is full and faithful.
\end{proof}

\subsection*{Extensivity}

Extensive sum types are semi-extensive finite sums for which the dependent types over a coproduct context are equivalently given by dependent types over the summand contexts.

\begin{definition}
  A display map category $(\mathbb B,\mathcal D)$ is \emph{extensive} if it has finite sum types and for each cospan of display maps $f : A \twoheadrightarrow I$, $g : B \twoheadrightarrow I$, the functor
    in \eqref{eq:ext}
  is an equivalence.
\end{definition}

Recall \cite[A1.4.4]{JOHNSTONE_elephantvols12} that a category $\mathbb C$ with finite coproducts \emph{has disjoint coproducts} if coproduct inclusions are monic and for any objects $A,B \in \mathbb C$, the following commuting square is a pullback.
\begin{displaymath}
  \begin{tikzpicture}
      \node (a) at (0,0) {$A$};
      \node (ab) at ({0.8*\arrow},0) {$A+B$};
      \node (z) at (0,{0.7*\arrow}) {$0$};
      \node (b) at ({0.8*\arrow},{0.7*\arrow}) {$B$};

      \path[->]
      (z) edge (b)
      (z) edge (a)
      (a) edge[right hook->] (ab)
      (b) edge[left hook->] (ab)
      ;
  \end{tikzpicture}
\end{displaymath}

\begin{proposition}
  Let $(\mathbb B,\mathcal D)$ be a well-rooted semi-extensive display map category.
  The following are equivalent.
  \begin{itemize}
  \item [(i)] $(\mathbb B,\mathcal D)$ is extensive.
  \item [(ii)] $\mathcal D$ is preserved by coproduct and in any diagram of the form \eqref{eq:10} where the vertical arrows are display maps and the bottom row is a coproduct diagram, if the top row is a coproduct diagram then both squares are pullbacks.
  \item [(iii)] In any diagram of the form \eqref{eq:10}
    where the\/ \emph{outer} vertical arrows are displays and the bottom row is a coproduct diagram, if the top row is also a coproduct diagram then the middle vertical arrow is a display map and both squares are pullbacks.
  \end{itemize}
  Moreover, if $(\mathbb B,\mathcal D)$ is extensive, then $\mathcal D$ contains the coproduct inclusions and coproducts are disjoint.
\end{proposition}

To tie extensive display map categories together with the usual notion of extensive category, note that a category $\mathbb C$ with pullbacks is extensive if and only if the class of all arrows in $\mathbb C$ is an extensive class of display maps.

\subsection*{Strong extensivity}

Strongly extensive finite sum types give us the expressivity we need for the $\Pi$-types in the polynomial model in Theorem \ref{theorem-polynomial}.

\begin{definition}
  A display map category $(\mathbb B,\mathcal D)$ is \emph{strongly extensive} if it is an extensive display map category and, for every object $X \in \mathbb B$, the category $\mathcal D/X$ is an extensive category.
\end{definition}

\begin{theorem}\label{theorem-strongly-extensive}
  Let $(\mathbb B,\mathcal D)$ be a well-rooted display map category.
  The following are equivalent.
  \begin{itemize}
  \item [(i)] $(\mathbb B,\mathcal D)$ is strongly extensive.
  \item [(ii)] $\mathbb B$ is an extensive category and $\mathcal D$ is preserved by copairing and coproduct.
  \end{itemize}
  If\/ $(\mathbb B,\mathcal D)$ has $\Sigma$-types, then we may include the following.
  \begin{itemize}
  \item [(iii)] $\mathbb B$ is an extensive category and $\mathcal D$ is preserved by copairing and contains all coproduct inclusions.
  \end{itemize}
\end{theorem}

\subsection*{Partial maps}

\begin{definition}
  Let $\mathbb C$ be a category with finite coproducts and $A$ and $B$ two objects of $\mathbb C$.
  A \emph{partial map} $A \rightharpoonup B$ is a tuple $(X,Y,i,j,f)$ where $X$ and $Y$ are objects in $\mathbb C$, $i : X \to A$ and $j : Y \to A$ are maps in $\mathbb C$ exhibiting $A$ as the coproduct of $X$ and $Y$, and $f : X \to B$ is a map in $\mathbb C$.
  Two partial maps $A \rightharpoonup B$ $(X,Y,i,j,f)$ and $(X',Y',i',j',f')$ are \emph{equivalent} when there exist isomorphisms $\theta_X : X \to X'$ and $\theta_Y : Y \to Y'$ such that $i' \circ \theta_X = i$, $j' \circ \theta_Y = j$, and $f' \circ \theta_X = f$.
\end{definition}

We always consider partial maps up to equivalence.
The following theorem is crucial in the construction of dependent products in the polynomial model in Theorem \ref{theorem-polynomial}.
It states that for any two types $A$ and $B$ there is a `partial function space' $R$, i.e.\ a type whose terms correspond to partial maps $A \rightharpoonup B$ rather than (total) morphisms $A \to B$.
Moreover this correspondence is given by substituting into a `generic partial function' $x : R \vdash (H,K,h,k,p) : A \rightharpoonup B$, and gives us a way in the type theory to talk about the domain and complement of the domain of a partial function.

\begin{theorem}\label{theorem-partial-maps}
  Let $(\mathbb B,\mathcal D)$ be a strongly extensive display map category with $\Sigma$- and\/ $\Pi$-types.
  Let $f : A \twoheadrightarrow I$ and $g : B \twoheadrightarrow I$ be two display maps.
  Then there exist a display map $r : R \twoheadrightarrow I$ and partial map $(h, k, i, j, p) : r^*(f) \rightharpoonup r^*(g)$ in $\mathcal D/R$, such that partial maps $f \rightharpoonup g$ in $\mathcal D/I$ correspond bijectively to sections $s$ of\/ $r$ via the operation sending a section $s$ to the coproduct decomposition $f \cong s^*r^*(f) \cong s^*(h) +_I s^*(k)$ and the map $s^*(p) : s^*(h) \to s^*r^*(g) \cong g$.
  Moreover, this bijection is natural in $g$ and stable under reindexing in $I$.
  Additionally, if\/ $(\mathbb B,\mathcal D)$ is well-rooted and if $f$ and $g$ are $u$-small for some universe $u$ closed under $\Sigma$-, $\Pi$-, and finite sum types, then $r$, $h$ and $k$ are all small maps.
\end{theorem}
\begin{proof}[Proof (sketch)]
  Take $r$ to be the fibred exponential $f \Rightarrow_I g +_I 1_I$.
  Now $h$ is given as the composite of the product projection $(f \Rightarrow_I g +_I 1_I) \times_I f \to f \Rightarrow_I g +_I 1_I$ with the pullback of the coproduct inclusion $g \hookrightarrow g +_I 1_I$ along the evaluation morphism $\mathsf{ev} : (f \Rightarrow_I g +_I 1_I) \times_I f \to g +_I 1_I$.
  We define $k$ using the other coproduct inclusion.
  The morphism $p$ is the evident one induced by the morphism $\mathsf{ev}^*(g) \to g$ from the pullback defining $h$.
  We omit the details of checking the bijection and naturality.

  For the last statement, we need only observe that coproduct inclusions are again small.
  Given types $x : X \twoheadrightarrow I$ and $y : Y \twoheadrightarrow I$, the inclusion $X \hookrightarrow X +_I Y$ is classified by $[\ulcorner 1_X \urcorner, \ulcorner 0_Y \urcorner] : X +_I Y \to \mathcal U$, where $\ulcorner 1_X \urcorner : X \to \mathcal U$ classifies the identity on $X$ and $\ulcorner 0_Y \urcorner : Y \to \mathcal U$ classifies the unique map $0_Y \twoheadrightarrow Y$.
\end{proof}

\section{The Dialectica or `polynomial' model}\label{sec:polynomial-models}

We are now ready to give the first of our Dialectica constructions, the \emph{polynomial model} introduced in \cite{VONGLEHN2014}.
The name, which we explain below, fits while we are considering the \emph{predicate-free} Dialectica construction.

\begin{definition}
  Let $(\mathbb B,\mathcal D)$ be a well-rooted display map category with finite sum types.
  Then the \emph{polynomial model} is the glued display map category $(\poly,\mathcal G)$, or just $\poly$, formed from the fibred display map category $(p : \poly \to \mathbb B,\mathcal E)$ over $(\mathbb B,\mathcal D)$ given as follows.
  The fibration $p$ is the opposite fibration to the codomain fibration $\mathcal D \to \mathbb B$.
  For each object $I \in \mathbb B$, the class $\mathcal E_I$ is the class of product projections in $\poly(I) = (\mathcal D/I)^\mathsf{op}$.
\end{definition}

Recall that $(\mathbb B,\mathcal D)$ having finite sum types means that $\mathsf{cod} : \mathcal D \to \mathbb B$ has fibred finite coproducts.
Hence $\poly \to \mathbb B$, being the opposite fibration, has fibred finite products.

\begin{lemma}
  The data $(p : \poly \to \mathbb B,\mathcal E)$ is indeed a fibred display map category with well-rooted fibres.
  Hence $\poly$ is indeed a well-rooted display map category.
\end{lemma}

Let us look more closely at $\poly$.
It is a version of the category of \emph{polynomials} or \emph{containers} \cite{AbbottAltenkirchGhani2003}, which has been shown to be cartesian closed \cite{AltenkirchLevyStaton2010}.
An object is simply a $\mathcal D$-map $f : A \twoheadrightarrow I$ in $\mathbb B$, (representing a Dialectica proposition $\exists\, i^I \forall a^A \top$ as in Section \ref{sec:dialectica-category}).
A morphism $(g : B \twoheadrightarrow J) \to (f : A \twoheadrightarrow I)$ consists of a pair $(h,\phi)$ making the diagram
\vspace{-0.3cm}
\begin{equation}\label{eq:6}
  \begin{tikzpicture}[baseline={0.4*\arrow cm}]
    \node (j) at (0,0) {$J$};
    \node (b) at (0,{0.8*\arrow}) {$B$};
    \node (i) at ({1.8*\arrow},0) {$I$};
    \node (a) at ({1.8*\arrow},{0.8*\arrow}) {$A$};
    \node[text depth=0] (ah) at ({0.8*\arrow},{0.8*\arrow}) {$A_h$};

    \path[->]
    (j) edge node[auto,swap] {$h$} (i)
    (ah) edge node[auto] {$h'$} (a)
    (ah) edge node[auto,swap] {$\phi$} (b)
    ;

    \path[->>]
    (a) edge node[auto] {$f$} (i)
    (ah) edge node[auto,outer sep=-1,pos=0.4] {$f'$} (j)
    (b) edge node[auto,swap] {$g$} (j)
    ;
  \end{tikzpicture}
\end{equation}
commute, where the inner square is a pullback.
A $\mathcal G$-map with codomain $f : A \twoheadrightarrow I$ is a morphism of the form
\begin{equation}\label{eq:7}
  \begin{tikzpicture}[baseline={0.4*\arrow cm}]
    \node (j) at (0,0) {$J$};
    \node[text depth=0] (b) at ({-0.3*\arrow},{0.8*\arrow}) {$A_h +_J X$};
    \node (i) at ({1.8*\arrow},0) {$I$};
    \node (a) at ({1.8*\arrow},{0.8*\arrow}) {$A$};
    \node[text depth=0] (ah) at ({0.8*\arrow},{0.8*\arrow}) {$A_h$};

    \path[->]
    (ah) edge node[auto,swap] {$\phi$} (b)
    ;

    \path[->>]
    (j) edge node[auto,swap] {$h$} (i)
    (ah) edge node[auto] {$h'$} (a)
    (a) edge node[auto] {$f$} (i)
    (ah) edge node[auto,outer sep=-1,pos=0.4] {$f'$} (j)
    (b) edge[shorten <=0.1cm] node[auto,swap] {$g$} (j)
    ;
  \end{tikzpicture}
\end{equation}
where $h \in \mathcal D$ and $g$ is the coproduct in $\mathcal D/J$ of $f'$ with some $x : X \twoheadrightarrow J$ and $\phi$ is the coproduct inclusion.

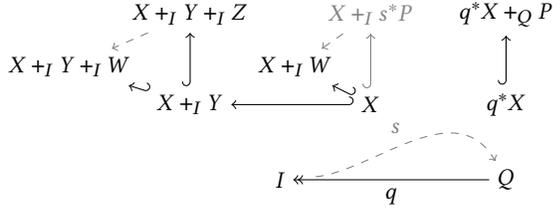
\begin{figure}
  \centering
  \begin{displaymath}
    \begin{tikzpicture}
      \node (i) at (0,0) {$I$};
      \node (x) at (1.2,1) {$X$};
      \node (xy) at (-1.2,1) {$X +_I Y$};
      \node (xyz) at (-1.2,2.2) {$X +_I Y +_I Z$};

      \node (q) at (3,0) {$Q$};
      \node (qx) at (3,1) {$q^* X$};
      \node (p) at (3,2.2) {$q^*X +_Q P$};

      \node[gray] (sp) at (1.2,2.2) {$X +_I s^*P$};

      \node (xw) at (0.2,1.5) {$X +_I W$};
      \node (xyw) at (-2.8,1.5) {$X +_I Y +_I W$};

      \path[->>]
      (q) edge node[auto] {$q$} (i)
      ;

      \path[dashed,->,gray]
      (i) edge[out=5,in=120,shorten <=0.3cm] node[auto] {$s$} (q)
      (sp) edge (xw)
      (xyz) edge (xyw)
      ;

      \path[gray]
      (x) edge[right hook->] (sp)
      ;

      \path
      (x) edge[left hook->] (xw)
      (xy) edge[left hook->] (xyw)
      (x) edge[left hook->] (xy)
      (xy) edge[right hook->] (xyz)
      (qx) edge[right hook->] (p)
      ;
    \end{tikzpicture}
  \end{displaymath}
  \caption{Quasifibred $\Pi$-types in $p : \poly \to \mathbb B$.}
  \label{fig:quasifibred-pi-poly}
\end{figure}

Let us give the main theorem.
\begin{theorem}\label{theorem-polynomial}
  Let\/ $(\mathbb B,\mathcal D)$ be a strongly extensive well-rooted display map category with $\Sigma$-, $\Pi$-, and identity types.
  Then $\poly$ is a strongly extensive well-rooted display map category with $\Sigma$-, $\Pi$-, and identity types.
  Moreover, if\/ $(\mathbb B,\mathcal D)$ has a universe closed under $\Sigma$-, $\Pi$-, identity, and finite sum types, then so does $\poly$.
\end{theorem}
\begin{proof}[Proof (sketch)]
  For $\Sigma$-types, we can simply apply Proposition \ref{proposition-fibred-sigma}, since product projections are closed under composition.

  For identity types we can apply Theorem \ref{theorem-fibred-identity-types} once we verify the left map condition (since left maps in the fibres are just the split monomorphisms, which are clearly preserved by reindexing).
  To do so, let $m : J \to I$ be a left map in $\mathbb B$, and let $x : X \twoheadrightarrow I$ and $y : Y \twoheadrightarrow I$ be two objects of $\poly(I)$, so that $X \hookrightarrow X +_I Y$ (in the opposite category) is the general form of a display map in $\poly(I)$, and we check that $m^*$ surjectively takes retractions (sections in the opposite category) of $X \hookrightarrow X +_I Y$ to retractions of $m^*(X) \hookrightarrow m^*(X) +_J m^*(Y)$.
  This amounts to, for any $h : m^*(Y) \to m^*(X)$, finding a dotted map in the following diagram.
  \begin{displaymath}
    \begin{tikzpicture}
      \node (i) at (3,0) {$I$};
      \node (j) at (0,0.4) {$J$};
      \node (y) at (2.5,0.7) {$Y$};
      \node (x) at (3.5,0.9) {$X$};
      \node (my) at (-0.5,1.1) {$m^*Y$};
      \node (mx) at (0.5,1.3) {$m^*X$};

      \path[->>]
      (mx) edge (j)
      (my) edge (j)
      (x) edge (i)
      (y) edge (i)
      ;
      \path[->]
      (my) edge (mx)
      (my) edge[crosses over] (y)
      (mx) edge (x)
      (j) edge node[auto,swap] {$m$} (i)
      ;
      \path[dotted,->]
      (y) edge[shorten >=0.1cm] (x)
      ;
    \end{tikzpicture}
  \end{displaymath}
  Since $(\mathbb B,\mathcal D)$ has $\Pi$-types, left maps are stable under pullback along $\mathcal D$-maps by \ref{lemma-identity-and-pi}, and hence we can use the left-lifting property of $m^*(Y) \to Y$ against $X \twoheadrightarrow I$.

  For $\Pi$-types, we observe that $\Sigma$-types in $(\mathbb B,\mathcal D)$ give us $\mathcal D$-sums in $\mathsf{cod} : \mathcal D \to \mathbb B$, and hence they give us $\mathcal D$-products in the opposite fibration.
  As $\mathcal D$-products are right adjoints, they preserve the fibrewise display maps, which are just product projections.
  Hence we can apply Theorem \ref{theorem-quasifibred-pi} once we verify that $p : \poly \to \mathbb B$ has quasifibred $\Pi$-types.

  Let $I \in \mathbb B$ and let $x : X \twoheadrightarrow I$, $y : Y \twoheadrightarrow I$, and $z : Z \twoheadrightarrow I$ be three objects in $\poly(I)$, so that a general composable pair of display maps is given by $X +_I Y +_I Z \hookleftarrow X +_I Y \hookleftarrow X$.
  Referring to Figure \ref{fig:quasifibred-pi-poly}, we need to find $q : Q \twoheadrightarrow I$ and $\pi : P \twoheadrightarrow Q$ such that $q$ together with $q^*X +_Q P \hookleftarrow q^*X$ form a quasifibred $\Pi$-type.
  We take an arbitrary $w : W \twoheadrightarrow I$, and note that the pullback in $\poly(I)$ of $X +_I W \hookleftarrow X$ along $X +_I Y \hookleftarrow X$ is given by pushout.
  Hence we must find $q$ and $\pi$ such that maps $X +_I Y +_I Z +_I \to X +_I Y +_I W$ over $X +_I Y$, i.e.\ maps $Z \to X +_I Y +_I W$, correspond to sections $s$ of $q$ together with a map $s^*P \to X +_I W$.
  We observe that the former kind of map is equivalently a partial map $f : Z \rightharpoonup Y$ together with a map $Z\backslash \mathop{\mathsf{dom}}(f) \to X +_I W$.
  Hence we can use Theorem \ref{theorem-partial-maps}: apply the theorem to $z : Z \twoheadrightarrow I$ and $y : Y \twoheadrightarrow I$ and take $q : Q \twoheadrightarrow I$ to be the $r : R \twoheadrightarrow I$ from the theorem, and $\pi : P \twoheadrightarrow Q$ to be the $k : K \twoheadrightarrow R$.

  For finite sums, the initial object of $\poly$ is $1_0 : 0 \to 0$ and the coproduct of $f : A \twoheadrightarrow I$ and $g : B \twoheadrightarrow J$ is just $f+g : A+B \twoheadrightarrow I+J$.
  To verify this, one uses Theorem \ref{theorem-strongly-extensive} and we omit the details, but note that we require the strongly extensive finite sums in $(\mathbb B,\mathcal D)$ to even get ordinary finite sums in $\poly$.

  Finally, we show that $\poly$ admits a universe, using Lemma \ref{lemma-quasifibred-universe} and constructing a quasifibred universe.
  Let $u : \widetilde{\mathcal U} \to \mathcal U$ be the universe in $\poly$.
  Take $\Omega = \mathcal U$, and take $v : \widetilde{\mathcal V} \twoheadrightarrow V$ in $\poly(\mathcal U)$ to be the map represented by $0 + \widetilde{\mathcal U} \hookleftarrow 0$.
  It is easy to check that this is indeed a quasifibred universe, and that the class of $v$-small maps is precisely the class of displays of the form \eqref{eq:7} where both $h : J \twoheadrightarrow I$ and $x : X \twoheadrightarrow J$ are $u$-small.
  We appeal to Theorem \ref{theorem-quasifibred-universe} (the verification of the hypotheses is easy) which still leaves us to check closure under finite sum types, but this is straightforward.
\end{proof}

Inputs to which Theorem \ref{theorem-polynomial} could be applied include in the simplest case the `extensional' models, i.e. locally cartesian closed and extensive categories where all maps are display maps. For example, the category of sets or more generally any topos, or the category of PERs. The result will be a (non-extensional) model satisfying the same conditions, so the construction can be iterated.
Other `intensional' examples include Voevodsky's simplicial model \cite{KLV_simplicial2012} and the cubical model of Bezem, Coquand, and Huber \cite{BEZEMCOQUANDHUBER2014}, which are easily checked to be strongly extensive.

\section{Biproducts and additive monads}\label{sec:biproducts}

Before constructing the Diller-Nahm model, we briefly review the notion of biproduct.
The following may be found in \cite{MACLANE_cwm1998}.

\begin{definition}
  A \emph{zero object} in a category $\mathbb C$ is an object $0$ which is both initial and terminal.
  For objects $A$, $B$ in a category $\mathbb C$ with zero object, the \emph{zero morphism} $A \to B$ is the unique morphism which factorizes through $0$.
\end{definition}

\begin{definition}
  In a category $\mathbb C$ with zero object, a \emph{biproduct} for objects $X$ and $Y$ consists of an object $X \oplus Y$ together with maps $\iota_X : X \to X \oplus Y$, $\iota_Y : Y \to X \oplus Y$, $\pi_X : X \oplus Y \to X$, and $\pi_Y : X \oplus Y \to Y$ satisfying the equations $\pi_X\iota_X = 1_X$, $\pi_Y\iota_Y = 1_Y$, $\pi_X\iota_Y = 0$, and $\pi_Y\iota_X = 0$, and moreover making $X \to X \oplus Y \leftarrow Y$ a coproduct diagram and $X \leftarrow X \oplus Y \to Y$ a product diagram.
\end{definition}

Let us see how biproducts can arise in the Kleisli category for a monad $M$ on some category $\mathbb C$ with finite products and finite coproducts.
Suppose that $M(0) \cong 1$, i.e.\ the initial object is sent to the terminal object.
Then for any objects $X$ and $Y$ it easy to construct a morphism $M(X + Y) \to MX \times MY$.
We take the following definition and theorem from \cite{COUMANSJACOBS_additivemonads2013}.

\begin{definition}[{\cite[Definition 7.16]{COUMANSJACOBS_additivemonads2013}}]
  A monad $M$ on a category $\mathbb C$ with finite products and finite coproducts is \emph{additive} if $M(0) \cong 1$ and the canonical map $M(X+Y) \to MX \times MY$ is an isomorphism.
\end{definition}

\begin{theorem}[{\cite[Theorem 7.20]{COUMANSJACOBS_additivemonads2013}}]
  Given an additive monad $M$ on a category $\mathbb C$ with finite products and finite coproducts, the Kleisli category $\mathbb C_M$ has biproducts, where the biproduct of $X$ and $Y$ is given by the coproduct $X + Y$ together with the $X$-inclusion $X \to X + Y \to M(X+Y)$ and the $X$-projection $X+Y \to M(X+Y) \to MX \times MY \to MX$, and similarly for $Y$.
\end{theorem}

\section{The Diller-Nahm model}\label{sec:diller-nahm-model}

The Diller-Nahm variant of the Dialectica interpretation \cite{DILLERNAHM1974} resolves a certain technical issue relating to the decidability of propositions.
For us, the difference is that now we assume our types are closed under passing to `finite multisets' $X \mapsto X^\bullet$ (the free commutative monoid) and a formula of the form \eqref{eq:1} instead reduces to
\begin{multline}\label{eq:9}
  \exists f^{U \to V}, F^{U \times Y \to X^\bullet} \forall u,y
  \\[-0.2cm]
  \Big(\big[\forall x \in F(u,y).\, A(u;x)\big] \to B(f(u);y)\Big),
\end{multline}
where we permit ourselves the use of quantification over a finite set inside the `quantifier-free' part of the formula.
The idea is that when $y$ renders $B(f(u);y)$ false, we may not have an effective procedure to decide which of a finite (multi)set $F(u,y)$ of possible counterexamples to $A(u;x)$ is indeed a counterexample, so we are required to check all of them.

The idea of a Diller-Nahm \emph{category} appears in \cite{DEPAIVA_thesis1991,HYLAND2002}.
This category has the same objects as $\mathsf{Dial}(\mathbb C)$, but now an arrow $(U,X,\alpha) \to (V,Y,\beta)$ is a realization of \eqref{eq:9}, i.e.\ it consists of arrows $f : U \to V$ and $F : U \times V \to M(X)$ in $\mathbb C$ together with a condition on $\alpha$ and $\beta$, where $M$ is an additive monad on $\mathbb C$ (as before, we simplify matters by ignoring the `$\alpha$' part of the objects).
Then this category $\mathsf{Dill}(\mathbb C)$ is already cartesian closed in this simply-typed case.
Our final result, Theorem \ref{theorem-diller-nahm}, is that the Diller-Nahm category generalizes to a model $\dill$ of dependent type theory which has $\Pi$-types, so in particular is still cartesian closed.
We will see how the biproducts in the Kleisli category for $M$ play a crucial role in constructing the $\Pi$-types in $\dill$.

For the Diller-Nahm model, we do not need to assume such strong extensivity properties as in the polynomial model, but we need to assume that we have an additive monad in each slice $\mathcal D/I$.
Moreover, we assume that it is a \emph{fibred monad}, meaning a monad $M$ on the total category $\mathcal D$ satisfying $\mathsf{cod} \circ M = \mathsf{cod}$ and whose unit and counit have $\mathsf{cod}$-vertical components.
Equivalently, this a family of monads $M_I$ on each category $\mathcal D/I$ which are preserved by the pullback functors $h^* : \mathcal D/J \to \mathcal D/I$ for each $h : I \to J$ in $\mathbb B$.
The monads represent the formation of a type of finite multisets.

\begin{definition}
  Let $(\mathbb B,\mathcal D)$ be a well-rooted display map category with finite sum types.
  Let $M$ be a fibred monad on the fibration $\mathsf{cod} : \mathcal D \to \mathbb B$, which is additive when restricted to a monad $M_I$ on each fibre category $\mathcal D/I$.
  Then the \emph{Diller-Nahm model} is the glued display map category $(\dill,\mathcal G)$, or just $\dill$, formed from the fibred display map category $(p : \dill \to \mathbb B)$ over $(\mathbb B,\mathcal D)$ given as follows.
  It is easy to check that the Kleisli category $(\mathcal D)_M$ is a fibred category over $\mathbb B$ with fibre category $(\mathcal D)_M(I)$ just the Kleisli category $(\mathcal D/I)_{M_I}$ of the original fibre category.
  The fibration $p$ is the opposite fibration to this fibred Kleisli category $(\mathcal D)_M \to \mathbb B$.
  For each object $I \in \mathbb B$, the class $\mathcal E_I$ is the class of product projections in $\dill(I) = ((\mathcal D/I)_{M_I})^\mathsf{op}$.
\end{definition}

The category $\dill$ has the same objects as $\poly$, but now a morphism $(g : B \twoheadrightarrow J) \to (f : A \twoheadrightarrow I)$ is a pair $(h,\phi)$ as in \eqref{eq:6} but $\phi$ represents a Kleisli arrow so is instead a morphism $f' \to M_J(g)$.
Likewise, a display map with codomain $f : A \twoheadrightarrow I$ is represented by a diagram of the form \eqref{eq:7}, but where $\phi$ is the coproduct (in fact, biproduct) inclusion in a Kleisli category.

\begin{theorem}\label{theorem-diller-nahm}
  Let $(\mathbb B,\mathcal D)$ be a well-rooted display map category with $\Sigma$-, $\Pi$-, identity, and finite sum types.
  Then $\dill$ is a well-rooted display map category with $\Sigma$-, $\Pi$-, and identity types.
  Moreover, if $(\mathbb B,\mathcal D)$ has a universe closed under $\Sigma$-, $\Pi$-, identity and finite sum types, then $\dill$ has a universe closed under $\Sigma$-, $\Pi$-, and identity types.
\end{theorem}
\begin{proof}[Proof (sketch)]
  The proof is very similar to that of Theorem \ref{theorem-polynomial}.
  The only interesting difference is in the construction of quasifibred $\Pi$-types, so we outline that here.
  Take the same setup as before and refer to Figure \ref{fig:quasifibred-pi-poly}, only now all the morphisms in the upper part of the diagram are Kleisli arrows, and the coproducts are moreover biproducts.
  We need to give $q : Q \twoheadrightarrow I$ and $\pi : P \twoheadrightarrow Q$ such that a Kleisli arrow $Z \to X \oplus_I Y \oplus_I W$ corresponds to sections $s$ of $q$ together with a Kleisli arrow $s^*P \to X \oplus_I W$, both Kleisli arrows being for the monad $M_I$ on $\mathcal D/I$.
  But by exploiting the biproducts, maps of the former kind correspond to pairs of Kleisli arrows $Z \to Y$ and $Z \to X \oplus_I W$.
  Hence we take $q : Q \twoheadrightarrow I$ to be the fibred exponential $(Z \Rightarrow_I M_IY) \twoheadrightarrow I$ and $\pi : P \twoheadrightarrow Q$ to be the pullback  $q^*(x +_I w)$ of $x +_I w : X +_I W \twoheadrightarrow I$ along $q$.
\end{proof}

\section{Conclusion}\label{sec:conclusion}

We have shown that the Dialectica construction generalizes from cartesian closed categories to categorical models of dependent type theory and presented two major examples which preserve $\Sigma$-, $\Pi$-, and identity types.
The dependently-typed setting even has the advantage over the simply-typed one that both constructions preserve cartesian closure.
Our proposed notion of strongly extensive finite sum types appears to be a fundamental one.
For reasons of space and clarity we have not presented either the Dialectica or Diller-Nahm models with a layer of predicates.
In fact these essentially rely only on the techniques we have developed in \S\ref{sec:gluing-construction}.
We have also omitted a third major example of a Dialectica construction, based on the error monad, which was considered for categories in \cite{BIERING_ccdc2008}, and which does require additional techniques to get a display map category with $\Pi$-types.
These constructions are considered in \cite{MOSS2017}, and we leave a presentation to future work.
We see our results here along with \cite{SHULMAN2014,UEMURA_fibredfibrationcategories2017} as the beginning of a \emph{model theory of dependent type theory}, wherein the gluing construction will be a sort of free completion (this is considered in \cite{VONGLEHN2014,MOSS2017}).

\begin{acks}
  We would like to thank Marcelo Fiore, Nicola Gambino, Martin Hyland, and Sam Staton for helpful discussions and comments.
  Sean Moss is currently supported by a Junior Research Fellowship at University College, Oxford, and previously by an EPSRC studentship at DPMMS, Cambridge.
  Tamara von Glehn is supported by a Junior Research Fellowship at Newnham College, Cambridge, and previously by a Cambridge International Scholarship from the Cambridge Overseas Trust.
\end{acks}

\bibliographystyle{ACM-Reference-Format}
\bibliography{dialectica-models.bib}

\end{document}